\documentclass[12pt]{article}
\usepackage[latin1]{inputenc}
\usepackage{amssymb,amsmath,amsfonts}
\newtheorem{theorem}{Theorem}
\newtheorem{proposition}[theorem]{Proposition}
\newtheorem{lemma}[theorem]{Lemma}

\newtheorem{corollary}[theorem]{Corollary}

\newcommand{\Ce}{\mathcal{C}}
\newcommand{\R}{\mathbb{R}}

\newcommand{\Sf}{\mathbb{S}}

\newcommand{\spa}{\mbox{span}}

\newcommand{\iis}{isometric immersions }

\def\D{{\cal D}}

\def\<{\langle}
\def\n{\nabla}

\def\>{\rangle}
\def\a{\alpha}

\def\bea{\begin{eqnarray*} }
\def\eea{\end{eqnarray*} }
\def\be{\begin{equation}}
\def\ee{\end{equation}}

\def\natural  #1{{\mathbb N^{#1}}}
\def\proof{\noindent{\it Proof: }}
\def\qed{\ifhmode\unskip\nobreak\fi\ifmmode\ifinner
\else\hskip5 pt \fi\fi\hbox{\hskip5 pt \vrule width4 pt
height6 pt  depth1.5 pt \hskip 1pt }}
\begin{document}

\title{Complete minimal submanifolds with nullity\\ in Euclidean space}
\author{M.\ Dajczer, Th.\ Kasioumis, A.\ Savas-Halilaj and Th.\ Vlachos}
\date{}
\maketitle
\renewcommand{\thefootnote}{\fnsymbol{footnote}} 
\renewcommand{\thefootnote}{\arabic{footnote}} 

\renewcommand{\thefootnote}{\fnsymbol{footnote}} 
\renewcommand{\thefootnote}{\arabic{footnote}}

\renewcommand{\thefootnote}{\fnsymbol{footnote}} 
\footnotetext{The third author would like to acknowledge financial support from the grant DFG SM 78/6-1.}  
\renewcommand{\thefootnote}{\arabic{footnote}} 

\begin{abstract} 
In this paper, we investigate minimal submanifolds in Euclidean 
space with positive index of relative nullity. Let $M^m$ be a complete 
Riemannian manifold and let $f\colon M^m\to\R^n$ be a minimal isometric 
immersion  with index of relative nullity at least $m-2$ at any point. 
We show that if the Omori-Yau maximum principle for the Laplacian holds 
on $M^m$, for instance, if the scalar curvature of $M^m$  does not decrease 
to $-\infty$ too fast or if the immersion $f$ is proper, then 
the submanifold must be a cylinder over a minimal surface.
\end{abstract}

\section{Introduction}

A frequent theme in submanifold theory is to find geometric conditions 
for an isometric immersion of a complete Riemannian manifold into Euclidean 
space $f\colon M^m\to\R^n$ with index of relative nullity $\nu\geq k>0$ at 
any point to be a $k$-cylinder. This means that the manifold $M^m$ splits as a 
Riemannian product $M^m=M^{m-k}\times\R^{k}$ and there is an isometric 
immersion $g\colon M^{m-k}\to\R^{n-k}$ such that
$f=g\times\rm{id}_{\R^{k}}$. 

The index of relative nullity introduced by Chern and Kuiper turned 
out to be a fundamental concept in the theory of isometric immersions.
At a point of $M^m$ the index  is just the dimension of 
the kernel of the second fundamental form of $f\colon M^m\to\R^n$ at that point. 
The kernels form an integrable distribution along any open subset where 
the index is constant and the images under $f$  of the leaves of the foliation  
are (part of) affine subspaces in the ambient space. Moreover, 
if $M^m$ is complete then the leaves are also complete along the open subset 
where the index reaches its minimum (cf.\ \cite{da}).
Thus, to conclude that $f$ is a cylinder one has to show that the images under 
$f$ of the leaves of relative nullity are parallel in the ambient space.

A fundamental result asserting that an isometric immersion 
$f\colon M^m\to\R^n$ with positive index of relative nullity  must be a 
$k$-cylinder is Hartman's theorem \cite{har} that requires the Ricci 
curvature of $M^m$ to be nonnegative; see also \cite{ma}.
A key ingredient for the proof of this result is the famous Cheeger-Gromoll 
splitting theorem used to conclude that the leaves of minimum relative nullity 
split intrinsically as a Riemannian factor.
Even for hypersurfaces, the same conclusion does not hold  if instead we 
assume that the Ricci curvature is nonpositive. 
Notice that the latter is always the case if $f$ is a minimal 
immersion. Counterexamples easy to construct are the complete irreducible 
ruled hypersurfaces of any dimension discussed in \cite[p.\ 409]{dg}.

Some of the many papers containing  characterizations of submanifolds as 
cylinders without the requirement of minimality are \cite{df,dg0,gf,har,ma,no,re}. 
When adding the condition of being  minimal we have 
\cite{abe,dr,fcz1,fcz2,gf,hsv,wz,yz}.

In this paper, we extend a result for hypersurfaces due to Savas-Halilaj 
\cite{sav} to the situation of arbitrary codimension.

\begin{theorem}\label{main}
Let $M^m$ be a complete Riemannian manifold and
$f\colon M^m\to\R^n$ be a minimal isometric immersion with index of relative 
nullity $\nu\geq m-2$ at any point of $M^m$. If the Omori-Yau maximum 
principle holds on $M^m$,  then $f$ is a cylinder over a minimal surface.
\end{theorem}

We recall that the \emph{Omori-Yau maximum principle} holds on $M^m$ if 
for any bounded from above function $\varphi\in C^\infty(M)$  
there exists a sequence of points $\{x_j\}_{j\in\natural{}}$
such that
$$
\lim\varphi(x_j)=\sup\varphi,\;\;\; \|\nabla\varphi\|(x_j)\leq 1/j
\;\;\;\text{and}\;\;\;\Delta\varphi (x_j)\leq 1/j
$$
for each $j\in\natural{}$. 

The category of complete Riemannian manifolds for which the  
principle is valid is quite large. For instance, it contains the 
manifolds with Ricci curvature bounded from below. 
It also contains the class of properly immersed submanifolds in a 
space form whose norm of the mean curvature vector is bounded 
(cf.\ \cite[Example 1.14]{pigola}).

\begin{corollary}\label{main2}
Let $M^m$ be a complete Riemannian manifold and $f\colon M^m\to\R^n$ be a 
minimal isometric immersion with index of relative nullity $\nu\geq m-2$ at any 
point of $M^m$. Assume that either the scalar curvature $\operatorname{scal}$
of $M^m$ satisfies $\operatorname{scal}\geq-c(d\log d)^2$ outside a compact 
set, where $c>0$ and  $d=d(\cdot,o)$ is the geodesic distance to a reference 
point $o\in M^m$, or that $f$ is proper. Then $f$ is a cylinder over a minimal 
surface.
\end{corollary}

Theorem \ref{main} is truly global in nature since there are plenty of (noncomplete) 
examples of minimal submanifolds of any dimension $m$ with constant index $\nu=m-2$ 
that are not part of a cylinder on any open subset. 
They can be all locally parametrically  described in terms of a certain 
class of elliptic surfaces; see Theorem $22$ in \cite{df}. In particular, there is
a Weierstrass type representation for these submanifolds when the manifold possesses 
a K\"ahler structure; see Theorem $27$ in \cite{df}. 
On the other hand, after the results of this paper what remains as a challenging 
open problem is the existence of a minimal complete noncylindrical submanifold  
$f\colon M^3\to\R^n$ with $\nu\geq 1$.

The main difficulty in the proof of Theorem \ref{main} arises from the fact that 
the index of relative nullity $\nu$ is allowed to vary. Consequently, one has to 
fully understand the structure of the set of points $\mathcal{A}\subset M^m$ where 
$f$ is totally geodesic in order to conclude that the  relative nullity foliation
on $M^m\smallsetminus\mathcal{A}$ extends smoothly to $\mathcal{A}$.

Recently Jost, Yang and Xin \cite{jost} proved various Bernstein type results for 
complete $m$-dimensional minimal graphical submanifolds in Euclidean space
with $\nu\ge m-2$. We observe that from a result in \cite{dg} 
it follows that the submanifolds considered in \cite[Theorem 1.1]{jost} are cylinders
over $3$-dimensional complete minimal submanifolds with $\nu\ge 1$.
Moreover, from Corollary \ref{main2} it follows that the submanifolds considered in 
\cite[Theorem 1.2]{jost} are just cylinders over complete minimal surfaces,
since entire graphs are  proper submanifolds. Thus, to prove a Bernstein theorem for such 
submanifolds is
equivalent to show a Bernstein theorem for entire minimal $2$-dimensional graphs 
in Euclidean space.

\section{Preliminaries}

In this first section, we recall some basic facts from the theory of \iis
that will be used in the proof of Theorem \ref{main}.

Let $M^m$ be a Riemannian manifold and $f\colon M^m\to\R^n$ be an isometric immersion.
As usual, often $M^m$ will be locally identified with its image. The \emph{relative nullity} 
subspace $\D(x)$ of $f$ at $x\in M^m$ is the  kernel of its second fundamental 
form $\alpha\colon TM\times TM\to N_fM$ with values in the normal bundle, 
that is,
$$
\D(x)=\{X\in T_xM:\alpha(X,Y)=0\;\;\text{for all}\;\;Y\in T_xM\}.
$$
Then, the dimension $\nu(x)$ of $\D(x)$ is called the \emph{index of relative nullity}  
of $f$ at $x\in M^m$. Let $U\subset M^m$ be an 
open subset where the index of relative nullity $\nu=s>0$ is  constant. 
It is a standard fact that the relative nullity distribution
$\D$ along $U$ is integrable, that the leaves of relative nullity are 
totally geodesic submanifolds of $M^m$ and that their images under $f$ are open 
subsets of affine subspaces in $\R^n$. 
The following is a well-known result in the theory of isometric 
immersions (cf.\ \cite[Theorem 5.3]{da}).

\begin{proposition}\label{comp} Let $\gamma\colon[0,b]\to M^m$ be a 
geodesic curve such that $\gamma([0,b))$ is contained in a leaf of relative nullity 
contained in $U$. Then also $\nu(\gamma(b))=s$.
\end{proposition}

The \emph{conullity space}  of $f$ at $x\in M^m$ is the orthogonal complement 
$\D^{\perp}(x)$ of $\D(x)$ in the tangent bundle $TM$. We write
$
X=X^v+X^h
$
according to the orthogonal splitting $TM=\D\oplus \D^{\perp}$
and denote ${\nabla}^{h}_XY = (\nabla_X Y )^h$.
The \textit{splitting tensor} $C\colon\D\times\D^{\perp}\to\D^{\perp}$ 
is given by
$$
C(T,X)=-{\nabla}^{h}_XT
$$
for any $T\in\D$ and $X\in\D^{\perp}$. 
The following differential equations for the tensor $C_T=C(T,\cdot)$
are well-known to hold (cf.\ \cite{da} or \cite{dg}):
\be\label{C1}
\n_S C_T=C_T C_S+C_{\n_ST}
\ee
and
\be\label{C2}
(\nabla^{h}_{X}C_T)Y-(\nabla^{h}_{Y}C_T)X
=C_{\nabla^{v}_{X}T}Y-C_{\nabla^{v}_{Y}T}X,
\ee
for any $S,T\in\Gamma(\D)$ and $X,Y\in\Gamma(\D^{\perp})$.
\vspace{1,5ex}

Finally, we have the following elementary result from the theory of 
submanifolds. 

\begin{proposition}\label{cylinder} Let $f\colon M^m\to\R^n$ be 
an isometric immersion with constant index of relative nullity $\nu=s>0$
and complete leaves of relative nullity. If the splitting tensor $C$
vanishes, then  $f$ is a $s$-cylinder.
\end{proposition}

\proof  That $C=0$ is equivalent to $\D$ being parallel in $M^m$. 
Consequently, the images via $f$ of the leaves of $\mathcal{D}$ are also
parallel in $\R^n$.\qed

\section{The proofs}

The possible structures of an isometric immersion $f\colon M^m\to\R^n$ when $M^m$ 
is complete and the index of relative nullity of $f$ satisfies $\nu\geq m-2$ at any 
point was completely described in \cite{dg}.  In particular, if $f$ is real 
analytic then it has to be either completely ruled or a cylinder
over a $3$-dimensional complete submanifold with $\nu\ge 1$.
In the case of minimal submanifols, it follows from Theorem $16$ in \cite{df} 
that we only have to consider the case of a nontrivial minimal
$f\colon M^3\to\R^n$ with $\nu\geq 1$ at any point of $M^3$. 

Let $U\subset M^3$ be an open subset where $\nu=1$ and the line bundle of 
relative nullity is trivial. 
Fix a smooth unit section $e$ spanning the relative nullity distribution  
along $U$ and let $J$ denote the unique, up to sign, almost complex 
structure acting on the conullity distribution $\D^{\perp}=\{e\}^\perp$. 
For simplicity, we set $\Ce=C_{e}$. Observe that our aim of proving 
Theorem \ref{main} will be achieved if we show that $\Ce$ is identically zero. 
The following lemma is of crucial importance.

\begin{lemma}\label{har}
There are harmonic functions $u,v\in C^{\infty}(U)$ such that
\be\label{C}
\Ce=v I-u J
\ee
where $I$ stands for the identity map on the conullity distribution.
\end{lemma}

\proof
We may assume that the immersion $f$ is substantial, that is, it does not reduce
codimension.
Let $A_{\xi}$ be the shape operator of $f$ with respect to the normal direction $\xi$, i.e.,
$$\langle A_{\xi}\,\cdot\,,\cdot\rangle=\<\alpha(\,\cdot\,,\cdot),\xi\>.$$
From the Codazzi equation for $A_{\xi}|_{\D^\perp}$ restricted to $\D^\perp$ we have that
$$
\nabla_{e}A_{\xi}|_{\D^\perp}= A_{\xi}|_{\D^\perp}\circ \Ce+A_{\nabla^{\perp}_e\xi}|_{\D^\perp}
$$
for any normal vector field $\xi\in N_fM$. Thus
$A_{\xi}|_{\D^\perp}\circ\Ce$ has to be symmetric, and hence
\be\label{cond1}
A_{\xi}|_{\D^\perp}\circ\Ce=\Ce^t\circ A_{\xi}|_{\D^\perp}.
\ee
On the other hand, the minimality condition is equivalent to 
\be\label{cond2}
A_{\xi}|_{\D^\perp}\circ J =J^t\circ A_{\xi}|_{\D^\perp}.
\ee

First we consider the hypersurface case $n=m+1$. 
Take a local orthonormal tangent frame $e_1,e_2,e_3$ that diagonalizes 
the shape operator of $f$ such that
$$Je_1=e_2\quad\text{and} \quad e_3=e$$
and let
$\xi$ be a unit normal along the hypersurface. Set
$$
u=\<\nabla_{e_2}e_1,e_3\>\;\;{\text {and}}\;\; v=\<\nabla_{e_1}e_1,e_3\>.
$$
From the Codazzi equation
$$
(\n_{e_i}A_{\xi})e_3=(\n_{e_3}A_{\xi})e_i,
$$
where $1\leq i\leq 2$, we have that $\<\nabla_{e_2}e_2,e_3\>=v$. Moreover, from
$$
\<(\n_{e_1}A_{\xi})e_2,e_3\>=\<(\n_{e_2}A_{\xi})e_1,e_3\>,
$$
we obtain that $\<\nabla_{e_1}e_2,e_3\>=-u$. Now we can readily see that (\ref{C}) holds true.

Now assume that $f$ is not an hypersurface. Consider the space
$$
N_1^f(x)=\spa\{\a(X,Y):\text{for all}\;\;X,Y\in T_xM\}.
$$
Notice that the dimension of $N_1^f(x)$ is at most two due to minimality.
Suppose that there is an open subset $V\subset M^3$ where $\dim N_1^f=1$.
A simple argument using the Codazzi equation \cite[Corollary 4.7]{da} shows 
that $N_1^f$ is parallel in the normal bundle along $V$, and thus the map $f\vert_V$ reduces 
codimension to an hypersurface. But due to real analyticity, the same would hold 
globally, and that is a contradiction.  Hence, there is an open dense subset $W$ of 
$M^3$ where $\dim N_1^f=2$. We conclude from (\ref{cond1}) and 
(\ref{cond2}) that $\Ce\in\spa\{I, J\}$ on $U\cap W$. By
continuity, we then get that $\Ce\in\spa\{I,J\}$ on $U$. Therefore, also in this case
there are functions $u,v\in C^{\infty}(U)$ such that (\ref{C}) holds.

It remains to show that  $u,v$ are harmonic. 
From (\ref{C1}) and (\ref{C2}) we have 
\be\label{c1}
\nabla^{h}_{e} \Ce=\Ce^2
\ee
and
\be\label{c2}
\big(\nabla^{h}_{X}\Ce\big)Y=\big(\nabla^{h}_{Y}\Ce\big)X
\ee
for any  $X,Y\in\Gamma(\D^{\perp})$.
For a local orthonormal tangent frame $e_1,e_2,e_3$ such that
$Je_1=e_2$ and $e_3=e$, it follows from (\ref{C}) that
\be\label{o1}
v=\<\nabla_{e_1} e_1,e_3\>=\<\nabla_{e_2}e_2,e_3\>
\ee
and
\be\label{o2}
u=-\<\nabla_{e_1}e_2,e_3\>=\<\nabla_{e_2}e_1,e_3\>.
\ee
It is easily seen that (\ref{c1}) is equivalent to 
\be\label{v}
e_3(v)=v^2-u^2\;\;\mbox{and}\;\; e_3(u)=2uv
\ee
whereas (\ref{c2}) to
\be\label{u}
e_1(u)=e_2(v)\;\;\mbox{and}\;\; e_2(u)=-e_1(v).
\ee
The Laplacian of $v$ is given by
\be\label{eq5}
\Delta v=\sum_{j=1}^3e_je_j(v)+\omega_{12}(e_2)e_1(v)-\omega_{12}(e_1)e_2(v)
-(\omega_{13}(e_1)+\omega_{23}(e_2))e_3(v)
\ee
where 
$$
\omega_{ij}(e_k)=\<\nabla_{e_k}e_i,e_j\>,
$$
where $1\leq i,j,k\leq 3$. Using (\ref{o2}) and (\ref{u}), we have that
\begin{align*}
e_1e_1(v)+e_2e_2(v)&=-e_1e_2(u)+e_2e_1(u)=[e_2,e_1](u)\\
&=\nabla_{e_2}e_1(u)-\nabla_{e_1}e_2(u)\\
&=\omega_{12}(e_1)e_1(u)+\omega_{12}(e_2)e_2(u)
+(\omega_{13}(e_2)-\omega_{23}(e_1))e_3(u)\\
&=\omega_{12}(e_1)e_2(v)-\omega_{12}(e_2)e_1(v)+2ue_3(u).
\end{align*}
Inserting the last equality into \eqref{eq5} and using (\ref{o1})
and \eqref{v} yields
$$
\Delta v= e_3e_3(v)+2u e_3(u)-2v e_3(v)=0.
$$
That also $u$ is harmonic is proved in a similar manner.\qed
\vspace{1ex}

Let us focus in the $3$-dimensional case, i.e., let $f\colon M^3\to\R^n$ 
be a minimal isometric immersion of a 
complete Riemannian manifold with index of relative nullity $\nu(x)\geq 1$ 
at any point $x\in M^3$, that is, the index is either $1$ or $3$. 
Let  $\mathcal{A}$ denote the set of totally geodesic points of $f$. 
 From Proposition~\ref{comp} the relative nullity foliation $\D$ 
is a line bundle on $M^3\smallsetminus\mathcal{A}$. 
Due to the real analyticity of the submanifold, the square of the norm of 
the second fundamental form is a real analytic function. It follows that
$\mathcal{A}$ is a real analytic set.  According to Lojasewicz's structure 
theorem \cite[Theorem 6.3.3]{kr} the set $\mathcal{A}$ locally decomposes as
$$
\mathcal{A}=\mathcal{V}^0\cup\mathcal{V}^1\cup \mathcal{V}^2\cup\mathcal{V}^3,
$$
where each $\mathcal{V}^d,\, 0\leq d\leq3$, is either empty or a 
disjoint finite union of $d$-dimensional real analytic subvarieties. 
A point $x_0 \in\mathcal{A}$ is called a \emph{regular point of dimension} $d$ 
if there is a neighborhood $\Omega$ of $x_0$ such that $\Omega\cap\mathcal{A}$  
is a $d$-dimensional real analytic submanifold of $\Omega$. 
If otherwise $x_0$ is said to be a \emph{singular} point.
The set of singular points is locally a finite union of submanifolds.

Our goal now is to  show that $\mathcal{A}=\mathcal{V}^1$, unless $f$ is
just an affine subspace in $\R^n$ in which case Theorem \ref{main} trivially holds. 
After excluding the latter trivial case, we have from the real analyticity of $f$
that $\mathcal{V}^3$ is empty.

\begin{lemma}\label{gp2}
The set $\mathcal{V}^2$ is empty.
\end{lemma}

\proof We only have to show is that there is no regular point in $\mathcal{V}^2$.
Suppose to the contrary that such a point do exist.
Let $\Omega \subset M^3$ be an open neighborhood of a smooth point $x_0\in\mathcal{V}^2$
such that $L^2=\Omega\cap\mathcal{A}$ is an embedded surface.
Let $e_1,e_2,e_3,\xi_1,...,\xi_{n-3}$ be an orthonormal frame  adopted to $M^3$
along $\Omega$ near $x_0$. The coefficients of the second fundamental form are
$$
h^a_{ij}=\<\a(e_i,e_j),\xi_a\>
$$
where from now on  $1\leq i,j,k\leq 3$ and $1\leq a,b\leq n-3$.

 The Gauss map $\gamma\colon M^3\to Gr(3,n)$ of $f$  as a map into 
the Grassmannian of oriented  $3$-dimensional subspaces in $\R^n$ 
is defined by $\gamma(x)=T_xM^3\subset\R^n$, up to parallel translation in $\R^n$ 
to the origin.  Regarding $Gr(3,n)$ as a submanifold in $\wedge^3\R^n$ 
via the map for the Pl\"{u}cker embedding, we have that
$\gamma=e_1\wedge e_2\wedge e_3$. Then
\be\label{gm}
\gamma_*e_i=\sum_{j,a}h^a_{ij}e_{ja}
\ee
where $ e_{ja}$ is obtained by replacing $e_j$ with $\xi_a$ in $e_1\wedge e_2\wedge e_3$.
Then
$$
\sum_i\<\gamma_*e_i,\gamma_*e_i\>=\sum_{i,j,a}(h^a_{ij})^2=\|\a\|^2
$$
where the inner product of two simple $3$-vectors in $\wedge^3\R^n$ is defined by
$$
\<a_1\wedge a_2\wedge a_3,b_1\wedge b_2\wedge b_3\>=\det\big(\<a_i,b_j\>\big).
$$
For a fixed simple $3$-vector $A=a_1\wedge a_2\wedge a_3$ let $w_A:M^3\to\R$
be the function defined by
$$
w_A=\<\gamma,A\>.
$$
Note that $w_A$ is a kind of height function.
Because the immersion $f$ is minimal, the function $w_A$ satisfies
$$
\Delta w_A=-\|\a\|^2w_A+\sum_{i,a\neq b,j\neq k}h^a_{ij}h^b_{ik}\<e_{ja,kb},A\>
$$
where $e_{ja,kb}$ is obtained by replacing $e_j$ with $\xi_a$ and $e_k$
with $\xi_b$ in $e_1\wedge e_2\wedge e_3$ (cf. \cite[p. 36]{x}).
Let $\varepsilon_1,\dots,\varepsilon_n$ be an orthonormal basis of $\R^n$ . 
The set 
$$
\left\{\varepsilon_{j_1}\wedge\varepsilon_{j_2}\wedge\varepsilon_{j_3}:
1\leq j_1<j_2<j_3 \leq n\right\}
$$
of $3$-vectors is an orthonormal basis  of  $\wedge^3\R^n$ by means of which  
identify $\wedge^3\R^n$ with $\R^{{ n \choose 3}}=\R^N$. 
Denoting by $\{A_J\}_{J\in\{1,\dots,N\}}$ the corresponding base in $\R^N$, we have
$$
\gamma=\sum_{J=1}^Nw_J A_J\;\;\mbox{where}\;\; w_J=\<\gamma, A_J\>.
$$
 From  $h^a_{ij}=\<\gamma_* e_i, e_{ja}\>$, we obtain
\be\label{h}
h^a_{ij}=\sum_{J}\< e_{ja},A_J\>e_i(w_J).
\ee
Moreover, for any $J\in\{1,\dots,N\}$, it holds
\be\label{d}
\Delta w_J=-\|\a\|^2w_J+\sum_{i,a\neq b,j\neq
k}h^a_{ij}h^b_{ik}\<e_{ja,kb},A_J\>.
\ee
Take a local chart $\phi\colon U\to\R^3$ of coordinates $x=(x_1,x_2,x_3)$ 
on an open subset $U$ of $\Omega$ and set
\be\label{e}
e_i=\sum_j\mu_{ij}{\partial_{x_j}}.
\ee
Setting $\theta_J=w_J\circ\phi^{-1}$, we obtain the map $\theta=\colon\phi(U)\subset\R^3\to \R^N$ 
given by
$$
\theta=\sum_J \theta_JA_J=(\theta_1,\dots,\theta_N).
$$
Note that $\theta=\gamma\circ\phi^{-1}$, i.e., $\theta$ is just the representation of the Gauss map
with respect to the above mentioned charts. From (\ref{h}) and (\ref{e}) we have 
\be\label{hh}
h^a_{ij}=\sum_{k,J}\mu_{ik}\< e_{ja},A_J\>(\theta_J)_{x_k}
\ee
and
\be\label{a}
\|\a\|^2
=\sum_{i,j,a}\Big(\sum_{k,J}\mu_{ik}\< e_{ja},A_J\> (\theta_J)_{x_k}\Big)^2.
\ee

The Laplacian of $M^3$ is given by
$$
\Delta=\frac{1}{\sqrt{g}} \sum_{i,j}{\partial_{x_i}}\Big(\sqrt{g}g^{ij}{\partial_{x_j}}\Big)
$$
where $g_{ij}$ are the components of the metric of $M^3$ and
$g=\det(g_{ij})$. Using (\ref{hh}) and (\ref{a}) we see  
that  (\ref{d}) is of the form
$$
\sum_{i,j}g^{ij} (\theta_J)_{x_ix_j}+C_J\big(x,\theta,\theta_{x_1},
\theta_{x_2},\theta_{x_3}\big)=0,
$$
where $C_J\colon\phi(U)\times\R^{4N}\to\R$ is given by
\bea
C_J(x,y,z_1,z_2,z_3)&=&\frac{1}{\sqrt g}
\sum_{i,j}(\sqrt{g}g^{ij})_{x_i}z_{jJ}+y_{J}
\sum_{i,j,a}\Big(\sum_{k,I}\mu_{ik}\<e_{ja},A_I\>z_{kI}\Big)^2\\
&&-\sum_{I,K}\sum_{i,l,m \atop a\neq b,j
\neq k}\mu_{il}\mu_{im}\<e_{ja,kb},A_J\>\<e_{ja},A_K\>
\<e_{kb},A_I\>z_{mI}z_{lK}
\eea
where $y=(y_1,\dots,y_N),z_i=(z_{i1},\dots,z_{iN})$, $i,m,l\in\{1,2,3\}$
and $I,J,K\in\{1,\dots,N\}$. 
Therefore, we have that the vector valued map $\theta=(\theta_1,\dots,\theta_N)$ satisfies the elliptic
equation
$$
\mathcal{L}\theta=\sum_{i,j}A_{ij}(x)\theta_{x_ix_j}+C\big(x,\theta,\theta_{x_1},\theta_{x_2},\theta_{x_3}\big)=0
$$
where $A_{ij}=g^{ij}I_N$, $I_N$ being the identity $N\times N$ matrix and
$C=(C_1,\dots,C_N)$.
Moreover, we have from (\ref{gm}) that  $\theta$ is constant on $\phi(L^2)$ 
and $\vec{n}(\theta)=0$ on $\phi(L^2)$ where  $\vec n$ is a unit normal field to the 
surface $\phi(L^2)$ in $\R^3$.

Consider the Cauchy problem $\mathcal{L}\theta=0$ with the following  initial conditions: 
$\theta$ is constant on $\phi(L^2)$ and $\vec n(\theta)=0$ on $\phi(L^2)$. 
According to the Cauchy-Kowalewsky  theorem (cf.\ \cite{t})
the problem has a unique solution if the surface $\phi(L^2)$ is
noncharacteristic. This latter is satisfied if  $Q(\vec n)\neq 0$,
where $Q$ is the characteristic form given by
$$
Q(\zeta)=\det (\Lambda(\zeta))
$$
where $\zeta=(\zeta_1,\zeta_2,\zeta_3)$ and
$$\Lambda(\zeta)=\sum_{i,j}g^{ij}\zeta_i\zeta_j I_N$$
is the symbol of the differential operator $\mathcal{L}$. That the surface $\phi(L^2)$  is 
noncharacteristic follows from
$$
Q(\zeta)= \Big(\sum_{i,j} g^{ij}\zeta_i\zeta_j\Big)^N.
$$
Because $C(x,y,0,0,0)=0$ the constant maps are solutions to the Cauchy problem. 
From the uniqueness part of the Cauchy-Kowalewsky theorem we conclude that the 
Gauss map $\gamma$ is constant on an open subset of $M^3$, and that is not
possible. \qed

\begin{lemma}\label{gp3}
The set $\mathcal{V}^0$ is empty.
\end{lemma}

\proof Suppose that $x_0\in\mathcal{V}^0$ and let $\Omega$ be an open 
neighborhood around  $x_0$ such that $\nu=1$ on  $\Omega\smallsetminus\{x_0\}$. 
Let $\{x_j\}_{j\in\natural{}}$ be a sequence in $\Omega\smallsetminus\{x_0\}$ 
converging to $x_0$. Let $e_j=e(x_j)\in T_{x_j} M$ be the 
sequence of unit vectors contained in the relative nullity distribution of $f$. 
By passing to a subsequence, if necessary, there is a unit vector 
$e_0\in T_{x_0} M$ such that $\lim e_j=e_0$. By continuity, the geodesic
tangent to $e_0$ at $x_0$ is a leaf of relative nullity outside $x_0$.
But this is a contradiction in view of Proposition \ref{comp}.\qed

\begin{lemma}\label{gp5}
The foliation $\mathcal{F}$ of the nullity distribution extends 
analytically over the regular points of $\mathcal{A}$.
\end{lemma}

\proof First observe that the relative nullity distribution  
extends continuously over the smooth points of $\mathcal{A}$.
In fact, by the previous lemmas  it remains to consider the case when $\Omega$ 
is an open subset of $M^3$ such that $\Omega\cap\mathcal{A}$  is a open segment
in a straight line in the ambient space. But in  this situation the result 
follows by a argument of continuity similar than in
the proof of Lemma \ref{gp3}.

Let $\Omega$ be an open subset of $M^3\smallsetminus\mathcal{A}$ 
and let $e_1,e_2,e_3$ be a local frame on $\Omega$ as in the proof
of Lemma \ref{har}. Consider the map  $F\colon\Omega\to\Sf^{n-1}$ into the unit 
sphere given by $F=f_*e_3$. A straightforward computation using 
(\ref{o1}), (\ref{o2}) and (\ref{u}) gives that its tension field
$$
\tau(F)=\sum_{j=1}^3\big(\overline\nabla_{F_*e_j}F_*e_j-F_*\nabla_{e_j}e_j\big)
$$
vanishes. Here $\overline{\nabla}$ denotes the
Levi-Civita connection of $\mathbb S^{n-1}$. Hence $F$ is a harmonic map.
Because $\mathcal{A}=\mathcal{V}^1$ its $2$-capacity ${\mathrm {cap}}_2(\mathcal{A})$
must be zero 
(cf.\ \cite[Theorem 3]{ev}). Since the map $F$ is continuous on $M^3$, it 
follows from a theorem of Meier \cite[Theorem $1$]{me}) that $F$ is of class
$C^2$ on $M^3$. But then $F$ is real analytic by a result due to Eells-Sampson 
\cite[Proposition  p.\ 117]{ee}.\qed

\begin{lemma}\label{gp6}
The set $\mathcal{A}$ has no singular points.
\end{lemma}

\proof According to Lemmas \ref{gp2} and \ref{gp3} the set
$\mathcal{A}$ only contains subvarieties of dimension one with possible
isolated singular points. Thus, by Lemma \ref{gp5}, the set of smooth points of 
$\mathcal{A}$ just contains segments of straight lines.
Hence, if there is a singular point in $\mathcal{A}$ it must be the intersection 
of such geodesic lines,  and that is clearly not possible. 
\vspace{1,5ex}\qed

The proof of our main result relies heavily on the following consequence
or the Omori-Yau maximum principle; see 
\cite[Theorem\ 28]{amr} or \cite[Lemma 4.1]{hsv2}.

\begin{lemma}\label{maxprinc2}
Let $M^m$ be a complete Riemannian manifold for which the Omori-Yau maximum 
principle holds. If $\varphi\in C^{\infty}(M)$  satisfies 
$\Delta\varphi\ge 2\varphi^2$ and $\varphi\geq 0$, then $\varphi = 0$.
\end{lemma}

\noindent\emph{Proof of Theorem \ref{main}:} 
Without loss of generality we may assume that $M^3$ is oriented by passing 
to the oriented double cover if necessary. 
It follows from Lemmas \ref{gp5} and \ref{gp6} that $J$ is
globally defined and that 
$\|\Ce\|^2=u^2+v^2$  is real analytic on $M^3$. 
From Lemma \ref{har} and \eqref{v} it follows that
\bea
\Delta (u^2+v^2)=2\|\nabla u\|^2 + 2\|\nabla v\|^2
\geq2(u^2+v^2)^2.
\eea
We deduce from Lemma \ref{maxprinc2} that $\Ce=0$, and by 
Proposition \ref{cylinder} this implies the desired splitting result.
\vspace{1,5ex}\qed

\noindent\emph{Proof of Corollary \ref{main2}:}  
The Omori-Yau maximum principle holds on $M^m$  
under the assumption on the scalar curvature (see \cite{aar} or 
\cite[Theorem 2.4]{amr}) or if the immersion $f$ is proper 
(see \cite[Theorem 2.5]{amr}).\qed

\noindent Marcos Dajczer\\
IMPA -- Estrada Dona Castorina, 110\\
Rio de Janeiro -- Brazil\\
e-mail: marcos@impa.br

\bigskip

\noindent Theodoros Kasioumis\\
University of Ioannina \\
Department of Mathematics\\
Ioannina--Greece\\
e-mail: theokasio@gmail.com

\bigskip

\noindent Andreas Savas-Halilaj\\
Leibniz Universit\"at Hannover \\
Institut f\"ur Differentialgeometrie\\
Welfengarten 1\\
30167 Hannover--Germany\\
e-mail: savasha@math.uni-hannover.de

\bigskip

\noindent Theodoros Vlachos\\
University of Ioannina \\
Department of Mathematics\\
Ioannina--Greece\\
e-mail: tvlachos@uoi.gr

\end{document}